\documentclass[11pt,dvipsnames,reqno]{amsart}
\usepackage[a4paper,left=3cm,right=3cm,top=3cm,bottom=3cm]{geometry}
\usepackage{amsmath,amsthm,amsfonts}
\newtheorem*{theorem*}{Theorem}
\usepackage{tikz}
\usepackage[T1]{fontenc}
\usepackage{microtype}
\usepackage{hyperref}
\definecolor{colorlinks}{RGB}{0, 24, 168}
\definecolor{colorcites}{RGB}{124, 10, 2}
\hypersetup{
    colorlinks=true,
    linkcolor=colorlinks,
    citecolor=colorcites,
    urlcolor=colorlinks,
    pdfborder={0 0 0}
}

\begin{document}


\title{The bunkbed conjecture on the complete graph}

\author{Peter van Hintum}
\address{Department of Pure Mathematics and Mathematical Statistics, University of Cambridge}
\email{pllvanhintum@maths.cam.ac.uk}

\author{Piet Lammers}
\address{Department of Pure Mathematics and Mathematical Statistics, University of Cambridge}
\email{p.g.lammers@statslab.cam.ac.uk}

\begin{abstract}
The bunkbed conjecture was first posed by Kasteleyn.
If $G=(V,E)$ is a finite graph and $H$ some subset of $V$, then the bunkbed of the pair $(G,H)$ is the graph $G\times\{1,2\}$ plus $|H|$ extra edges to connect for every $v\in H$ the vertices $(v,1)$ and $(v,2)$.
The conjecture asserts that $(v,1)$ is more likely to connect with $(w,1)$ than with $(w,2)$ in the independent bond percolation model for any $v,w\in V$.
This is intuitive because $(v,1)$ is in some sense closer to $(w,1)$ than it is to $(w,2)$.
The conjecture has however resisted several attempts of proof.
This paper settles the conjecture in the case of a constant percolation parameter and $G$ the complete graph.
\end{abstract}


\maketitle

\section{Introduction and main result}
The bunkbed conjecture is an intuitive statement in percolation theory.
In rough terms the conjecture asserts that -- in a specific setting
and in a specific sense -- two vertices of a graph
are more likely to remain connected after randomly removing some edges
if the graph distance between the vertices is smaller.
The conjecture is appealing because it is intuitive yet difficult to prove.
In this paper we prove the conjecture for the case that the underlying graph
is symmetrical.
The conjecture was first posed by Kasteleyn (in 1985),
as was remarked by Van den Berg and Kahn \cite{van2001correlation}.
Before stating the conjecture, we introduce the notion
of the bunkbed of a graph and we introduce the percolation model.

First, let $G=(V,E)$ be a finite graph, and let $H$ be a subset of $V$.
The bunkbed of the pair $(G,H)$, or $\operatorname{BB}(G,H)$,
is the graph $G\times\{1,2\}$ plus $|H|$
extra edges to connect for every $v\in H$
the vertices $(v,1)$ and $(v,2)$.
For any vertex $v\in V$, write $v^-:=(v,1)$ and $v^+:=(v,2)$.
Any vertex of $\operatorname{BB}(G,H)$ is of the form $v^-$ or $v^+$.
Equivalently, if $e\in E$, then write $e^\pm$ for the two corresponding edges
in the bunkbed graph.

Now introduce the bond percolation model for the bunkbed graph.
Pick a percolation parameter $p\in [0,1]$.
In the percolation model, every edge of the form $e^\pm$
is declared open with probability $p$ and closed
with probability $1-p$, independently of the other edges.
The edges of the form $\{v^-,v^+\}$ are always declared open.
Write $\mathbb P_p$ for the measure corresponding to the states of the edges.
For $v,w\in V$, write $v\sim w$ if $\{v,w\}$ is an open edge,
and write $v\leftrightarrow w$ if $v$ and $w$ are joined by an open path.
Furthermore, if $v\in V$ and $W\subset V$, then write $v\sim W$ if there is
a vertex $w\in W$ with $v\sim w$.
See \cite{grimmett2018probability} for a more elaborate introduction into the percolation
model.

The general bunkbed conjecture asserts that
$ \mathbb P_p(v^-\leftrightarrow w^-)$ is larger than or equal to $\mathbb P_p(v^-\leftrightarrow w^+)$,
for any $v,w\in V$.
This is precisely how the conjecture is described in \cite{haggstrom2003probability}.
We prove the conjecture in the case that $G$ is the complete graph.
Write $K_{n}$ for the complete graph on the vertex set $[n]:=\{1,2,...,n\}$.

\begin{theorem*}
  Pick $n\in\mathbb N$ and $H\subset [n]$.
  Consider independent bond percolation on $\operatorname{BB}(K_n,H)$
  with parameter $p\in [0,1]$ for the edges
  of the form $e^\pm$, and with the edges of the form $\{v^-,v^+\}$ always open.
  Then for any pair of vertices $v,w\in [n]$
  we have
  \begin{equation}
    \label{theorem}
    \mathbb P_p(v^-\leftrightarrow w^-)\geq\mathbb P_p(v^-\leftrightarrow w^+).
  \end{equation}
\end{theorem*}

De Buyer proved the theorem for the special case $p=\frac12$ in \cite{de2016proof},
which he later extended to $p\geq\frac 12$ in
\cite{de2018proof}\footnote{The proof in the current paper was presented at a seminar
in the Instituto de Matem\'atica Pura e Aplicada
in Rio de Janeiro on 8 February 2018.
It seems thus reasonable to assume that the case $p>\frac12$ was settled
simultaneously and independently by de Buyer and by the authors of the current paper.}.
The proof presented here draws on a different method.
The conjecture has been proved for any $p$ for wheel graphs and some small other graphs
by Leander \cite{leander2009sjalvstandiga}
and for outerplanar graphs by Linusson \cite{linusson2011percolation}.
Both \cite{leander2009sjalvstandiga} and \cite{linusson2011percolation} use the method of minimal counterexamples.
It has been proved that the connection probability of two vertices of a graph is the same
in the percolation model with parameter $p=\frac{1}{2}$
as it is in the model in which every edge is assigned a direction uniformly at random
\cite{karp1990transitive,linusson2011percolation,mcdiarmid1981general}.
A statement similar to the bunkbed conjecture has been studied on bunkbed graphs.
Bollob{\'a}s and Brightwell considered a continuous time random walk on a bunkbed graph,
such that the jump rate to any neighbour of the current state
is one \cite{bollobas1997random}.
They conjectured that for every $t>0$, this random walk
started at $v^-$ is more likely to have hit $w^-$ than $w^+$ before time $t$.
H{\"a}ggstr{\"o}m proved this conjecture in \cite{haggstrom1998conjecture}.

\section{Proof of the theorem}
\begin{proof}[Proof of the theorem]
We prove the theorem for $n+1$ instead of $n$
for notational convenience (the conjecture is trivial for $n=1$).
  It will be assumed that $w=n+1$, without loss of generality.
If $w\in H$ then $w^-\sim w^+$ and the two events in (\ref{theorem})
are the same.
If $v=w$ then the left side of (\ref{theorem}) equals one.
Therefore we only need to consider the case that $w\not \in H$ and $v\neq w$.
If $v\in H$ then both sides of (\ref{theorem}) are equal (by symmetry of the bunkbed),
and if $v\not\in H$, then both sides of (\ref{theorem})
do not depend on the actual choice of $v\in [n]\setminus H$ (by symmetry of the complete graph).
Therefore it is sufficient to prove the inequality for $v$ chosen uniformly at random
in the set $[n]$, independently of the percolation.
By choosing $v$ uniformly at random in $[n]$, we make optimal use of the symmetry of
the graph $K_{n+1}$.

Now write $(V,E):=\operatorname{BB}(K_{n+1},H)$ and note that
$V=[n+1]\times \{1,2\}=([n]\times\{1,2\} )\cup \{w^-,w^+\}$.
Write $O$ for the open subgraph of $\operatorname{BB}(K_{n+1},H)$
induced by the set $[n]\times\{1,2\}$.
This means that the vertex set of $O$ is $[n]\times\{1,2\}=V\setminus \{w^-,w^+\}$,
and that every edge $e\in E$ is an edge of $O$
if and only if its endpoints are in $[n]\times\{1,2\}$
and if $e$ is open in the percolation measure $\mathbb P_p$.
The edge set of $O$ is thus random in the measure $\mathbb P_p$.
Moreover, $O$ determines the configuration of all edges
incident to neither $w^-$ nor $w^+$, and the configuration of the edges
incident to either $w^-$ or $w^+$ and the value of $v$ are independent of $O$.
Write $c$ for the partition of $O$ into connected components,
and label these $c=\{c_1,...,c_k\}$
(where $k$ is the number of connected components, also random).

In order to calculate the difference between the two probabilities in (\ref{theorem}), we define the events
\begin{align*}
  A&:=\{w^-\leftrightarrow v^- \not\leftrightarrow w^+\}
=\dot\cup_i\left( \{v^-\in c_i\}\cap\{w^-\sim c_i\}\cap\{w^+\not\sim c_i\}\cap\{\not\exists j\neq i,\,w^-\sim c_j\sim w^+\}\right),\\
B&:=\{w^+\leftrightarrow v^- \not\leftrightarrow w^-\}
=\dot\cup_i \left(\{v^-\in c_i\}\cap\{w^+\sim c_i\}\cap\{w^-\not\sim c_i\}\cap\{\not\exists j\neq i,\,w^-\sim c_j\sim w^+\}\right).
\end{align*}
In each of these two equations the four events within the
disjoint union are, conditional on $O$ and for fixed $i$,
mutually independent.
This is because the last three events (in each of the two lines) depend on the states of different edges,
and because (for the first event) the value of $v$ is chosen independently of the percolation.
Write $\tilde{\mathbb P}$
for the measure $\mathbb P_p$ conditioned on $O$.
Now
\begin{alignat*}{3}
&\mathbb P_p(A|O)&&=\textstyle\sum_i
\tilde{\mathbb P}(v^-\in c_i)
\tilde{\mathbb P}(w^- \sim c_i)
\tilde{\mathbb P}(w^+ \not\sim c_i)
&&\tilde{\mathbb P}(\not\exists j\neq i,\,w^-\sim c_j\sim w^+),
\\
&\mathbb P_p (B|O)&&=\textstyle\sum_i
\tilde{\mathbb P}(v^-\in c_i)
\tilde{\mathbb P}(w^+ \sim c_i)
\tilde{\mathbb P}(w^- \not\sim c_i)
&&\tilde{\mathbb P}(\not\exists j\neq i,\,w^-\sim c_j\sim w^+),
\\
&\mathbb P_p(A|O)-\mathbb P_p(B|O)
&&=\textstyle\sum_i
\tilde{\mathbb P}(v^-\in c_i)
\left(
\tilde{\mathbb P}(w^+ \not\sim c_i)-
\tilde{\mathbb P}(w^- \not\sim c_i)
\right)
&&\tilde{\mathbb P}(\not\exists j\neq i,\,w^-\sim c_j\sim w^+).
\end{alignat*}
The difference between the two sides of (\ref{theorem})
is $\mathbb P_p(A)-\mathbb P_p(B)$, which equals
the expectation of the final line of the display over $O$.
The probabilities in (\ref{theorem}) are invariant under simultaneously replacing $v^-$ by $v^+$
and interchanging $w^-$ and $w^+$.
Taking the average over the original expression and the permuted one gives
\begin{align*}
&\mathbb P_p(v^-\leftrightarrow w^-)-\mathbb P_p(v^-\leftrightarrow w^+)=\mathbb P_pA-\mathbb P_pB
=\mathbb E_p \left(\mathbb P_p(A|O)-\mathbb P_p(B|O)\right)\\
&\quad=\frac{1}{2}
\mathbb E_p\textstyle\sum_i
\begin{array}{r}
  \Big(\quad
  \tilde{\mathbb P}(v^-\in c_i)
  \left(
  \tilde{\mathbb P}(w^+ \not\sim c_i)-
  \tilde{\mathbb P}(w^- \not\sim c_i)
  \right)
  \tilde{\mathbb P}(\not\exists j\neq i,\,w^-\sim c_j\sim w^+)
  \phantom{\Big)}
  \\
  +\:
  \tilde{\mathbb P}(v^+\in  c_i)
  \left(
  \tilde{\mathbb P}(w^- \not\sim  c_i)-
  \tilde{\mathbb P}(w^+ \not\sim  c_i)
  \right)
  \tilde{\mathbb P}(\not\exists j\neq i,\,w^+\sim   c_j\sim w^-)\Big)
\end{array}
\\&\quad=
\frac12
\mathbb E_p\textstyle\sum_i
\left(\tilde{\mathbb P}(v^-\in c_i)
-\tilde{\mathbb P}(v^+\in c_i)
\right)
\left(
\tilde{\mathbb P}(w^+ \not\sim c_i)-
\tilde{\mathbb P}(w^- \not\sim c_i)
\right)
\tilde{\mathbb P}(\not\exists j\neq i,\,w^-\sim c_j\sim w^+).
\end{align*}
We claim that the two differences in the final sum always have the same sign,
so that the product is always nonnegative.
Write $c_i^-$ and $c_i^+$ for the number of vertices in $c_i$ of the form
$u^-$ and $u^+$ respectively, so that, for example, $\sum_ic_i^-=\sum_ic_i^+=n$.
We explicitly calculate that (writing $q=1-p$)
$$\left(\tilde{\mathbb P}(v^-\in c_i)
-\tilde{\mathbb P}(v^+\in c_i)
\right)
\left(
\tilde{\mathbb P}(w^+ \not\sim c_i)-
\tilde{\mathbb P}(w^- \not\sim c_i)
\right)=\frac1n(c_i^--c_i^+)\left(q^{c_i^+}-q^{c_i^-}\right)\geq 0,$$
where the final inequality is due to $(a-b)\left(q^b-q^a\right)\geq 0$
for any $a,b\in\mathbb Z_{\geq 0}$ and $q\in [0,1]$.
\end{proof}

\section*{Acknowledgement}
The authors thank B\'ela Bollob\'as
for suggesting the problem.

The authors thank the referee for carefully reading the manuscript.

The first author was supported by the Cambridge Trust
and the UK Engineering and Physical Sciences Research Council.
The second author was supported by the Department of
Pure Mathematics and Mathematical Statistics, University of Cambridge and the UK
Engineering and Physical Sciences Research Council grant EP/L016516/1.

\bibliographystyle{amsplain}
\bibliography{ms}

\providecommand{\bysame}{\leavevmode\hbox to3em{\hrulefill}\thinspace}
\providecommand{\MR}{\relax\ifhmode\unskip\space\fi MR }
\providecommand{\MRhref}[2]{%
  \href{http://www.ams.org/mathscinet-getitem?mr=#1}{#2}
}
\providecommand{\href}[2]{#2}
\begin{thebibliography}{10}

\bibitem{van2001correlation}
Jacob van~den Berg and Jeff Kahn, \emph{A correlation inequality for connection
  events in percolation}, The Annals of Probability \textbf{29} (2001), no.~1,
  123--126.

\bibitem{bollobas1997random}
B{\'e}la Bollob{\'a}s and Graham Brightwell, \emph{Random walks and electrical
  resistances in products of graphs}, Discrete applied mathematics \textbf{73}
  (1997), no.~1, 69--79.

\bibitem{de2016proof}
Paul~de Buyer, \emph{A proof of the bunkbed conjecture for the complete graph
  at $p=\frac12$}, arXiv:1604.08439v1 (2016).

\bibitem{de2018proof}
\bysame, \emph{A proof of the bunkbed conjecture on the complete graph for
  $p\geq\frac12$}, arXiv:1802.04694v1 (2018).

\bibitem{grimmett2018probability}
Geoffrey Grimmett, \emph{Probability on graphs}, Cambridge University Press,
  2018.

\bibitem{haggstrom1998conjecture}
Olle H{\"a}ggstr{\"o}m, \emph{On a conjecture of {B}ollob{\'a}s and
  {B}rightwell concerning random walks on product graphs}, Combinatorics,
  Probability and Computing \textbf{7} (1998), no.~4, 397--401.

\bibitem{haggstrom2003probability}
\bysame, \emph{Probability on bunkbed graphs}, Proceedings of FPSAC, vol.~3,
  2003.

\bibitem{karp1990transitive}
Richard Karp, \emph{The transitive closure of a random digraph}, Random
  Structures \& Algorithms \textbf{1} (1990), no.~1, 73--93.

\bibitem{leander2009sjalvstandiga}
Madeleine Leander, \emph{Sj{\"a}lvst{\"a}ndiga arbeten i matematik, on the
  bunkbed conjecture},  (2009).

\bibitem{linusson2011percolation}
Svante Linusson, \emph{On percolation and the bunkbed conjecture},
  Combinatorics, Probability and Computing \textbf{20} (2011), no.~1, 103--117.

\bibitem{mcdiarmid1981general}
Colin McDiarmid, \emph{General percolation and random graphs}, Advances in
  Applied Probability \textbf{13} (1981), no.~1, 40--60.

\end{thebibliography}
\end{document}